\documentclass[a4paper,reqno]{amsart}
\numberwithin{equation}{section}

\def\Fbbd{\mathbb{F}}

\def\Zbbd{\mathbb{Z}}
\def\Acal{\mathcal{A}}

\def\Fcal{\mathcal{F}}
\def\Gcal{\mathcal{G}}
\def\Hcal{\mathcal{H}}

\def\kbf{\mathbf{k}}

\def\qbf{\mathbf{q}}
\def\wbf{\mathbf{w}}

\def\zbf{\mathbf{z}}

\def\phibf{\boldsymbol{\varphi}}
\def\xibf{\boldsymbol{\xi}}
\def\etabf{\boldsymbol{\eta}}

\def\ffrak{\mathfrak{f}}
\def\gfrak{\mathfrak{g}}
\def\a{\alpha}
\def\b{\beta}
\def\g{\gamma}
\def\d{\delta}

\def\la{\lambda}
\def\La{\Lambda}
\def\s{\sigma}
\def\phi{\varphi}

\def\abs#1{\left|#1\right|}

\def\Ref#1{(\ref{#1})}
\def\tilde{\widetilde}
\def\endproof{\hfill\rule{2mm}{2mm}}

\def\at{\tilde{a}}
\def\bt{\tilde{b}}
\def\ct{\tilde{c}}
\def\vt{\tilde{v}}
\def\ut{\tilde{u}}
\def\wt{\tilde{w}}
\def\st{\tilde{s}}
\def\tt{\tilde{t}}
\def\ft{\tilde{f}}
\def\gt{\tilde{g}}
\def\htilde{\tilde{h}}
\def\kt{\tilde{k}}
\def\mt{\tilde{m}}
\def\nt{\tilde{n}}
\def\pt{\tilde{p}}
\def\rt{\tilde{r}}
\def\xt{\tilde{x}}
\def\yt{\tilde{y}}
\def\zt{\tilde{z}}
\def\itilde{{\tilde{\iota}}}
\def\St{\tilde{S}}
\def\Tt{\tilde{T}}
\def\endproof{\hfill\rule{2mm}{2mm}}
\def\?{(?)\marginpar{|?}}

\newtheorem{theo}{Theorem}[section]
\newtheorem{prop}{Proposition}[section]

\def\beq{\begin{equation}}
\def\eeq{\end{equation}}
\def\be{\begin{equation*}}
\def\ee{\end{equation*}}
\begin{document}
\title{Eigenproblem for Jacobi matrices: hypergeometric series solution}
\author{Vadim B. Kuznetsov}
  \address{Department of Applied Mathematics,
          University of Leeds,
          Leeds LS2 9JT, UK}
  \email{V.B.Kuznetsov@leeds.ac.uk}
\author{Evgeny K. Sklyanin}
\address{Department of Mathematics, University of York,
York YO10 5DD, UK}
\email{eks2@york.ac.uk }
%
%
%
\begin{abstract}
We study the perturbative power-series expansions
of the eigenvalues and eigenvectors of a general tridiagonal (Jacobi)
matrix of dimension $d$. The (small) expansion parameters are the entries
of the two diagonals of length $d-1$ sandwiching the principal
diagonal which gives the unperturbed spectrum.

The solution is found explicitly
in terms of multivariable (Horn-type) hypergeometric series
in $3d-5$ variables in the generic case. To derive the result,
we first rewrite the spectral problem for the Jacobi matrix
as an equivalent system of algebraic equations
which are then solved by the application
of the multivariable Lagrange inversion formula.
The corresponding Jacobi determinant is calculated explicitly.
Explicit formulae are also found for any monomial composed
of eigenvector's components.

\end{abstract}
\maketitle
\tableofcontents

\section{Introduction}\label{intro}
\noindent
The problem of solving algebraic equations by series expansions
has a long history. In case of a single polynomial equation 
a major cornerstone is Birkeland's paper \cite{Bir1927}
where a solution is given through an application of
the Lagrange inversion formula \cite{AAR}.
Later on, an alternative approach was suggested by Mayr in \cite{Mayr1937},
who derived the relevant hypergeometric series as solutions
of some PDE's satisfied by the zeros as functions of the coefficients of
the polynomial.
A modern interpretation in terms of Gel'fand-Kapranov-Zelevinsky
hypergeometric functions \cite{GKZ1989,GKZ1990,GKZ1994} can be found in \cite{Sturm2000}. 
The latter approach can also be applied to general systems of algebraic equations.

The main goal of this paper is to derive a complete power-series solution of the
spectral problem for a finite ($d\times d$) Jacobi matrix $M$. 
We consider the off-diagonal matrix elements
to be small, so that the whole problem looks like a perturbation
of a diagonal matrix.
After fixing a normalization of the eigenvector
the problem is reduced to solving a system of $d$ {\it quadratic}
equations for $d$ unknowns. Then we transform it into
an equivalent {\it larger} system of special (Lagrange-form)
$3d-5$ {\it cubic} equations, which are then inverted by the application
of the multivariable Lagrange inversion formula.
The expansion of an arbitrary monomial of the components of the eigenvector
is given explicitly
in terms of multivariable (Horn-type) hypergeometric series
of $3d-5$ variables. In the special case of an eigenvalue 
growing from a corner matrix element
the number of expansion variables drops to $2d-3$.

Consider a tridiagonal (Jacobi) matrix of order $d$
\beq
  M=
   \begin{pmatrix}
       \a_1 & \b_1 \\
     \g_1 & \a_2 & \b_2 \\
    & \hdotsfor{3} \\
    && \g_{k-1} & \a_k & \b_k \\
    &&& \hdotsfor{3} \\
    &&&& \g_{d-1} & \a_d
  \end{pmatrix}
\eeq
and the corresponding eigenproblem  $MV=\La V$ for the eigenvector
$V$ and the eigenvalue $\La$.

Assume that the off-diagonal elements $\b_k$, $\g_k$ are the small parameters
of the power expansion and that $\a_k$ are distinct.
In the zeroth approximation $\b_k=\g_k=0$ the matrix $M$
is diagonal, its eigenvalues and eigenvectors being $\a_k$ and $V^{(k)}$, respectively,
where the components $V^{(k)}_j$ of $V^{(k)}$ can be chosen as $V^{(k)}_j=\d_{jk}$.
By continuity argument, for the small values of $\b_k$, $\g_k$
the eigenvalues $\La_k$, $k=1,\ldots,d$, are distinct and can be numbered
in such a way that
\beq
\La_k=\a_k+\text{higher order terms}.
\label{branch_lambda}
\eeq
We also choose to normalise the eigenvector $V^{(k)}$
\beq
    MV^{(k)}=\La_k V^{(k)}
\label{Mv=lambdav}
\eeq
by the condition
\beq
V^{(k)}_k=1
\label{normv}
\eeq
for its $k^{\text{th}}$ component. Therefore, the remaining components must
vanish in the zeroth approximation:
\beq
    V^{(k)}_j=0+\text{higher order terms}, \qquad j\neq k.
\label{branch_v}
\eeq

The eigenvalue problem \Ref{Mv=lambdav} together with
the normalization condition \Ref{normv} produces 
a system of $d$ algebraic (quadratic) equations
for the eigenvalue $\La_k$ and the free components $V^{(k)}_j$, $j\neq k$
of the eigenvector $V^{(k)}$
defining them as algebraic functions of the parameters $\a$, $\b$, $\g$.
The conditions \Ref{branch_lambda} and \Ref{branch_v} fix uniquely
the branches of the multivalued algebraic functions for small values of
$\b$ and $\g$.

The problem we solve in the present paper is to find an effective way to construct
explicit expressions for the coefficients of the power series expansions
for the eigenvalues $\La_k$ and the components of the eigenvectors $V^{(k)}_j$.
The available literature on solving systems of algebraic equations
by multivariate hypergeometric series (see a review in \cite{Sturm2000})
focuses mainly on solving the generic algebraic system.
We are not aware of any detailed analysis of
the particular systems arising from the Jacobi matrix spectral problem.
The importance of the latter problem for numerous applications 
in mathematical physics and, in particular, in the theory of quantum
integrable systems
has been the main motivation of our study.
Rather then using more modern approaches,
in the present paper we follow the original idea of Birkeland \cite{Bir1927} and use
a variant of the Lagrange inversion formula. 

We shall use the following variant of
the multivariable Lagrange inversion theorem \cite{Good1960,Ges1978}.
Let boldface letters
denote vectors $\xibf\equiv(\xi_1,\ldots,\xi_D)$, 
multi-indices $\boldsymbol{q}\equiv(q_1,\ldots,q_D)\in\Zbbd^D$, and monomials
$\boldsymbol{\xi^q}\equiv \xi_1^{q_1}\ldots \xi_D^{q_D}$. 
The inequality $\boldsymbol{q}\geq\boldsymbol{p}$ is understood component-wise:
$q_i\geq p_i$, $\forall i$.
Let $[\boldsymbol{\xi^q}]h(\xibf)$ denote
the coefficient at $\boldsymbol{\xi^q}$ in the power series $h(\xibf)$.

\begin{theo}\label{th:lagrange}
Let $\xibf=(\xi_1,\ldots,\xi_D)$, $\etabf=(\eta_1,\ldots,\eta_D)$.
Let $\phi_i(\etabf)$, $i=1,\ldots,D$, be
formal power series in $\etabf$ such that $\phi_i(\mathbf{0})\neq0\,\forall i$.
Then the system of $D$ equations
\beq
   \xi_i=\frac{\eta_i}{\phi_i(\etabf)}, \qquad i=1,\ldots,D,
\label{zPhi}
\eeq
defines uniquely $\eta_i(\xibf)$,
$i=1,\ldots,D$, as formal power series in $\xibf$.

In addition, let $\chi(\etabf)$ be a multiple Laurent series,
that is $\chi(\etabf)\boldsymbol{\eta^p}$ is a power series for some
$\boldsymbol{p}\geq0$. Then the Laurent series expansion for
$\chi\bigl(\boldsymbol{\eta(\xi)}\bigr)$ 
\beq
    \chi\bigl(\boldsymbol{\eta(\xi)}\bigr)=\sum_{\boldsymbol{q}\geq\boldsymbol{-p}}
    X_{\boldsymbol{q}}\,\boldsymbol{\xi^q}
\label{laurent_series}
\eeq
is given by the formula
\beq
  X_{\boldsymbol{q}}\equiv[\boldsymbol{\xi^q}]\chi\bigl(\boldsymbol{\eta(\xi)}\bigr)
   =[\boldsymbol{\eta^q}]\chi(\etabf)\boldsymbol{\phi^q(\eta)}J(\etabf),
\label{lagrange_expansion}
\eeq
where $J$ is the Jacobian
\beq
   J= \det\left(\mathbf{1}-\frac{\partial \log \phibf}{\partial \log\etabf}\right)
   =\det\left|\d_{jk}-\frac{\eta_k}{\phi_j}\,\frac{\partial \phi_j(\etabf)}{\partial \eta_k}\right|.
\label{defJ}
\eeq
\end{theo}

The analytic version of the Lagrange theorem
\cite{Good1960} guarantees that 
if $\phi_i(\etabf)$ and $\chi(\etabf)\boldsymbol{\eta^p}$ are analytic at $0$ then
the Laurent series \Ref{laurent_series} 
converges in a punctured neighbourhood around $\xibf=\boldsymbol{0}$.

When applying the Lagrange formula \Ref{lagrange_expansion},
the major complication comes from the
Jacobian $J$, which may be difficult to compute.
Fortunately, for our particular problem the Jacobian can be calculated explicitly,
in a relatively compact form.

The paper is organised as follows. In section \ref{lagrange} we
write down the set of quadratic equations defining the eigenvalue $\La$
and the components $V_j$ of the eigenvector $V$, transform them into the form which
is convenient for study, identify the combinations of the small parameters
which serve as the expansion variables, and rewrite the equations again
in the form which allows us to apply Lagrange's inversion formula.
In section \ref{jacobian} we calculate an important component of
Lagrange's formula: the Jacobian $J$.

In section \ref{hypergeometric} we put together all the ingredients of
the Lagrange formula and produce explicit expressions for all $\Lambda_k$'s and
$V_j^{(k)}$'s as finite sums of power series. The number of terms in the sum
equals to the number of terms in the Jacobian $J$. All the power series are
particular cases of a single universal Horn-type multivariable hypergeometric
series which we denote $\Phi$. In generic situation, the series $\Phi$ depends on 
$3d-5$ expansion variables and $d-1$ integer parameters.
In section \ref{corner} we describe the simplification of our
results  for the special case $k=1$ (or $k=d$) when the eigenvalue
$\La$ stems from a corner of the matrix $M$. 
In section \ref{examples} we examine a few low-dimensional examples
illustrating the general results.
The last section contains
a discussion of possible applications and extensions of our result.

\section{Lagrange equations}\label{lagrange}
\noindent
{}From now on we shall concentrate on studying a single
eigenvalue $\La_k$ and the corresponding eigenvector
$V^{(k)}$, for some fixed value of the index $k$ \Ref{branch_lambda}.
We shall change our notation accordingly, to simplify the calculations.
Set $r\equiv d-k$ and $\rt\equiv k-1$, so that 
\beq
   d=r+\rt+1.
\label{def_d}
\eeq

Let $\La_k=\a_k+\la$, $V^{(k)}_j=\vt_{k-j}$ and $\a_j=\at_{k-j}$
for $j=1,\ldots,k$, and
$V^{(k)}_j=v_{j-k}$ and $\a_j=a_{j-k}$
for $j=k,\ldots,d$, so that $V^{(k)}_k=v_0=\vt_0$
and $\a_k=a_0=\at_0$.
Respectively, let $\b_j=\bt_{k-j-1}$ and $\a_j=\at_{k-j-1}$
for $j=1,\ldots,k-1$,
and $\b_j=\b_{j-k}$ and $\a_j=\a_{j-k}$
for $j=k,\ldots,d-1$.
Without loss of generality we can set $\a_k=a_0=\at_0=0$.
The normalisation condition \Ref{normv} implies that $V^{(k)}_k=v_0=\vt_0=1$.
Since $\a_j$ are assumed distinct, we have $a_i\neq0$ and $\at_i\neq0$
for $i\neq0$.

As a result, the eigenvalue problem \Ref{Mv=lambdav} takes the following form:
\beq
  \left(\begin{array}{ccccccccccccc}
   \at_{\rt} & \bt_{\rt-1} \\
   & \hdotsfor{3} \\
   & \ct_i & \at_i & \bt_{i-1} \\
   && \hdotsfor{3} \\
   &&& \ct_1 & \at_1 & \bt_0 \\
   &&&& \ct_0 & 0 & b_0  \\
   &&&&&  c_0 & a_1 & b_1 \\
   &&&&&& \hdotsfor{3} \\
   &&&&&&& c_{i-1} & a_i & b_i \\
   &&&&&&&& \hdotsfor{3} \\
   &&&&&&&&& c_{r-1} & a_r
  \end{array}\right)
  \begin{pmatrix}
   \vt_{\rt} \\ \ldots \\ \vt_i \\ \ldots \\ \vt_1 \\ 1 \\ v_1 \\ \hdots \\ v_i \\ \hdots \\ v_r
  \end{pmatrix}
  =\lambda
  \begin{pmatrix}
   \vt_{\rt} \\ \ldots \\ \vt_i \\ \ldots \\ \vt_1 \\ 1 \\ v_1 \\ \hdots \\ v_i \\ \hdots \\ v_r
  \end{pmatrix},
\label{jacobi}
\eeq
The numbers $\rt,r=0,1,\ldots$ measure the distances of the selected diagonal element
from the corners of the matrix. The cases $\rt,r=0,1$ are slightly special,
we shall comment on them in due course, all other cases with $\rt,r\ge 2$
are generic. 

Expressing $\la$ from the `central' (`zeroth') row as
\beq
\la=\ct_0\vt_1+b_0v_1
\label{lav1}
\eeq
and substituting it into the remaining rows we get the set of $r+\rt=d-1$, see 
\Ref{def_d}, quadratic equations
\begin{subequations}
\label{eq_vvti}
\begin{alignat}{2}
    c_{i-1}v_{i-1}+a_iv_i+b_iv_{i+1}&=\ct_0\vt_1v_i+b_0v_1 v_i, &\qquad
      i&=1,\ldots,r,
\label{eq_vi} \\
   \ct_i\vt_{i+1}+\at_i\vt_i+\bt_{i-1}\vt_{i-1}&=\ct_0\vt_1\vt_i+b_0v_1\vt_i
      &\qquad i&=1,\ldots,\rt,
\label{eq_vti}
\end{alignat}
\end{subequations}
where we assume that $v_0\equiv\vt_0=1$, $b_r=0$, $\ct_{\rt}=0$. Formula
\Ref{lav1} eliminates the eigenvalue $\la$ by expressing it in terms 
of the two eigenvector's components: $v_1$ and $\vt_1$. From now on, the variables $v_i$ 
and $\vt_i$ will be our only $d-1$ unknowns.

Note that equations \Ref{eq_vvti}
are invariant with respect to the `involution'
(rotation by 180 degrees around the central element)
\beq
r \leftrightarrow\rt, \quad
v_i\leftrightarrow\vt_i, \quad
a_i\leftrightarrow\at_i, \quad
b_i\leftrightarrow\ct_i, \quad
c_i\leftrightarrow\bt_i
\label{tilde_sym}
\eeq
which we refer to as the {\em tilde-symmetry}. 

The next step is to rescale the variables $v_i$ and $\vt_i$
in order to identify the convenient combinations of the expansion parameters.
Setting $b_0=\ldots=b_{r-1}=0$ and $\ct_0=\ldots=\ct_{\rt-1}=0$
and denoting the corresponding values of
$v_i$ and $\vt_i$ as, respectively, $v_i^0$ and $\vt_i^0$ we get the equations
\begin{subequations}
\begin{alignat}{2}
    c_{i-1}v_{i-1}^0+a_iv_i^0&=0 &\qquad \Longrightarrow \qquad
    v_i^0&=-\frac{c_{i-1}}{a_i}\,v_{i-1}^0,\\
   \at_i\vt_i^0+\bt_{i-1}\vt_{i-1}^0&=0 &\qquad \Longrightarrow \qquad
    \vt_i&=-\frac{\bt_{i-1}}{\at_i}\,\vt_{i-1}^0,
\end{alignat}
\end{subequations}
which can be solved recursively yielding
\begin{subequations}
\begin{alignat}{2}
    v_i^0&=(-1)^i\,\frac{c_0\ldots c_{i-1}}{a_1\ldots a_i},
   &\qquad i&=1,\ldots,r,\\
    \vt_i^0&=(-1)^i\,\frac{\bt_0\ldots \bt_{i-1}}{\at_1\ldots \at_i},
   &\qquad i&=1,\ldots,\rt.
\end{alignat}
\end{subequations}

Let us rescale
$v_i$ and $\vt_i$ by the formulae
\begin{subequations}\label{vfromu}
\begin{alignat}{2}
    v_i&=v_i^0u_i=(-1)^i\,\frac{c_0\ldots c_{i-1}}{a_1\ldots a_i}\,u_i,
   &\qquad i&=1,\ldots,r,\\
    \vt_i&=\vt_i^0\ut_i=(-1)^i\,\frac{\bt_0\ldots \bt_{i-1}}{\at_1\ldots \at_i}\,\ut_i,
   &\qquad i&=1,\ldots,\rt
\end{alignat}
\end{subequations}
and rewrite equations \Ref{eq_vvti} in terms of $u_i$, $\ut_i$:
\begin{subequations}
\label{eq_uuti}
\begin{alignat}{2}
    u_i&=u_{i-1}-\frac{b_0c_0}{a_1a_i}\,u_1u_i
         -\frac{\bt_0\ct_0}{\at_1a_i}\,\ut_1u_i
         +\frac{b_ic_i}{a_ia_{i+1}}\,u_{i+1},
    &\qquad  i&=1,\ldots,r,
\label{eq_ui}\\
    \ut_i&=\ut_{i-1}-\frac{\bt_0\ct_0}{\at_1\at_i}\,\ut_1\ut_i
           -\frac{b_0c_0}{a_1\at_1}\,u_1\ut_i
    +\frac{\bt_i\ct_i}{\at_i\at_{i+1}}\,\ut_{i+1},
    &\qquad  i&=1,\ldots,\rt,
\label{eq_uti}
\end{alignat}
\end{subequations}
where we assume $u_0=\ut_0=1$, $b_r=\ct_r=0$.

Equations \Ref{eq_uuti} contain the parameters
$\boldsymbol{a,b,c,\at,\bt,\ct}$ in specific combinations, which are
convenient to use as the expansion parameters.
Note that $b_i$ and $c_i$ enter only through the product $b_ic_i$
(same for $\bt_i\ct_i$).
Introduce the variables $x_i$, $y_i$, $z_i$
\begin{subequations}
\label{def_xyz}
\begin{alignat}{2}
  x_i&=\frac{b_0c_0}{a_1a_{i+1}}, &\qquad i&=0,\ldots,r-1, \\
  y_i&=\frac{\bt_0\ct_0}{\at_1a_{i+1}}, &\qquad i&=0,\ldots,r-1,\\
  z_i&=\frac{b_{i+1}c_{i+1}}{a_{i+1}a_{i+2}}, &\qquad i&=0,\ldots,r-2,
\end{alignat}
\end{subequations}
(altogether $r+r+(r-1)=3r-1$ variables) as well as their tilde-analogs
\begin{subequations}
\label{def_xyzt}
\begin{alignat}{2}
  \xt_i&=\frac{\bt_0\ct_0}{\at_1\at_{i+1}}, &\qquad i&=0,\ldots,\rt-1, \\
  \yt_i&=\frac{b_0c_0}{a_1\at_{i+1}}, &\qquad i&=0,\ldots,\rt-1,\\
  \zt_i&=\frac{\bt_{i+1}\ct_{i+1}}{\at_{i+1}\at_{i+2}}, &\qquad i&=0,\ldots,\rt-2,
\end{alignat}
\end{subequations}
(altogether $\rt+\rt+(\rt-1)=3\rt-1$ variables).
The total number of variables is then
$D=3r+3\rt-2=3d-5$. 

Equations \Ref{eq_uuti} are simplified now to the form
\begin{subequations}
\label{eq_uuti_xyz}
\begin{alignat}{2}
    u_i&=u_{i-1}-x_{i-1}u_1u_i-y_{i-1}\ut_1u_i+z_{i-1}u_{i+1},
    &\qquad  i&=1,\ldots,r,
\label{eq_ui_xyz}\\
    \ut_i&=\ut_{i-1}-\xt_{i-1}\ut_1\ut_i-\yt_{i-1}u_1\ut_i+\zt_{i-1}\ut_{i+1},
    &\qquad  i&=1,\ldots,\rt,
\label{eq_uti_xyz}
\end{alignat}
\end{subequations}
where we assume $u_0=\ut_0=1$, $z_{r-1}=\zt_{\rt-1}=0$.
Since $b_i,\bt_i,c_i,\ct_i$ are small parameters, so are $x_i,y_i,z_i,\xt_i,\yt_i,\zt_i$.
In the zeroth approximation we have a set of binomial equations
$u_i^0=u_{i-1}^0$, $\ut_i^0=\ut_{i-1}^0$.
The iteration of equations \Ref{eq_uuti_xyz} with the initial values
$u_i^0=\ut_i^0=1$ produces formal power series expansions of $u_i$ and $\ut_i$
in the variables $(\boldsymbol{x,y,z,\xt,\yt,\zt})\equiv\boldsymbol{\xi}$.

The tilde-symmetry \Ref{tilde_sym} for equations \Ref{eq_uuti_xyz}
takes the form
\beq
r \leftrightarrow\rt, \quad
u_i\leftrightarrow\ut_i, \quad
x_i\leftrightarrow\xt_i, \quad
y_i\leftrightarrow\yt_i, \quad
z_i\leftrightarrow\zt_i.
\label{tilde_sym__u}
\eeq

{\em Important remark:} the variables $\boldsymbol{\xi}$
are, generally speaking, not independent. Indeed, the $3d-5$ variables
are bound by the relations
\beq
    \frac{x_i}{y_i}=\frac{\yt_j}{\xt_j}=\frac{b_0c_0\at_1}{\bt_0\ct_0a_1},
    \qquad i=0,\ldots,r-1, \quad
    j=0,\ldots,\rt-1,    
\eeq
altogether $r+\rt-1=d-2$ independent relations, which leaves only $2d-3$
independent variables. Only in the special case $r=0$ (or $\rt=0$), 
see end of this section, the variables $\boldsymbol{\xi}$ become independent.
Nevertheless, when solving equations \Ref{eq_uuti_xyz} by power series,
it is convenient to treat $\boldsymbol{\xi}$ as the set of $3d-5$ independent variables.
The expressions \Ref{def_xyz}--\Ref{def_xyzt} for $\boldsymbol{\xi}$
in terms of $(\boldsymbol{a,b,c,\at,\bt,\ct})$ can be substituted then
into the resulting expansions at the very final stage.

In order to find explicitly the coefficients of the power series for
$u_i$ and $\ut_i$ we shall use the Lagrange inversion formula.
In the notation of Theorem \ref{th:lagrange}, the expansion variables vector
$\xibf$ is composed of the 6 sets of variables $\boldsymbol{\xi=(x,y,z,\xt,\yt,\zt)}$.
To meet the theorem's premises we have to introduce also the vector $\etabf$
composed, respectively, of the 6 matching sets of expandable quantities
$\boldsymbol{\eta=(s,t,w,\st,\tt,\wt)}$:
\begin{subequations}
\begin{alignat}{4}
    s_i&\sim x_i &\qquad i&=0,\ldots,r-1;&\qquad
    \st_j&\sim \xt_j &\qquad j&=0,\ldots,\rt-1, \\
    t_i&\sim y_i &\qquad i&=0,\ldots,r-1;&\qquad
    \tt_j&\sim \yt_j &\qquad j&=0,\ldots,\rt-1, \\
    w_i&\sim z_i &\qquad i&=0,\ldots,r-2;&\qquad
    \wt_i&\sim \zt_i &\qquad i&=0,\ldots,\rt-2.
\end{alignat}
\end{subequations}

Note that the variables $(\boldsymbol{u,\ut})$ cannot be used directly as $\etabf$
because they have nonzero limits for $\xibf\rightarrow \boldsymbol{0}$,
therefore they have to be completed by small factors.
Besides, the number of the variables $(\boldsymbol{u,\ut})$ is only $r+\rt=d-1$,
so we need $2d-4$ extra variables to match the number $D=3d-5$
of expansion variables. Actually, it is more convenient to
split the number $D$ differently:
$D\equiv 3d-5=(2d-2)+(d-3)$, with the first set of $2d-2$ variables
for the larger system coming from only two variables, $u_1$ and $\ut_1$.
Define $s_i$, $t_i$, $\st_i$, $\tt_i$  as
\begin{subequations}
\label{def_sst0}
\begin{align}
    s_i&=x_iu_1, &t_i& =y_i \ut_1,&& i=0,\ldots,r-1,
\label{def_s0}\\
    \st_i&=\xt_i\ut_1, &\tt_i& =\yt_i u_1,&& i=0,\ldots,\rt-1.
\label{def_st0}
\end{align}
\end{subequations}
The second set of $d-3$ remaining variables $\boldsymbol{w,\wt}$ 
is defined as
\begin{subequations}
\label{def_wwti}
\begin{alignat}{2}
    w_i=z_i\frac{u_{i+2}}{u_i}, \qquad i=0,\ldots,r-2,
\label{def_wi}\\
    \wt_i=\zt_i\frac{\ut_{i+2}}{\ut_i}, \qquad i=0,\ldots,\rt-2.
\label{def_wti}
\end{alignat}
\end{subequations}

It might seem easier to define $w_i$ to be proportional to $u_i$ rather than to the ratio
of two $u$'s with the step 2. However, a little experimenting shows that our choice
of $w_i$ leeds to simpler, factorised expressions for the functions $\phi_i(\boldsymbol{\eta})$
in \Ref{zPhi}.

{}From $u_0=\ut_0=1$ it follows that
\begin{subequations}\label{exxx}
\begin{align}
    w_0=z_0u_2,
\label{def_w0}\\
    \wt_0=\zt_0\ut_2.
\label{def_wt0}
\end{align}
\end{subequations}
Solving equations \Ref{def_wwti}--\Ref{exxx} recursively we obtain the
expressions for $u_2,\ldots,u_r$ 
in terms of $w_j$ and for $\ut_2,\ldots,\ut_{\rt}$ in terms of $\wt_j$:
\begin{subequations}\label{ufromw}
\begin{alignat}{2}
  u_{2j}&= \frac{w_{2j-2}\ldots w_0}{z_{2j-2}\ldots z_0}, &\qquad
  u_{2j+1}&= \frac{w_{2j-1}\ldots w_1s_0}{z_{2j-1}\ldots z_1x_0},\\
  \ut_{2j}&= \frac{\wt_{2j-2}\ldots w_0}{\zt_{2j-2}\ldots \zt_0}, &\qquad
  \ut_{2j+1}&= \frac{\wt_{2j-1}\ldots \wt_1\st_0}{\zt_{2j-1}\ldots \zt_1\xt_0}.
\end{alignat}
\end{subequations}

\begin{prop} \label{pprop}
The set of equations \Ref{eq_uuti_xyz}, \Ref{def_sst0} and \Ref{def_wwti}
is equivalent to the set of Lagrange-type \Ref{zPhi} equations
\begin{subequations}
\begin{alignat}{5}
   x_i&=\frac{s_i}{f_i}, &\qquad
   y_i&=\frac{t_i}{g_i}, &\qquad i&=0,\ldots,r-1, &\qquad
   z_j&=\frac{w_j}{h_j}, &\qquad j&=0,\ldots,r-2,
\label{def_fgh}\\
   \xt_i&=\frac{\st_i}{\ft_i}, &\qquad
   \yt_i&=\frac{\tt_i}{\gt_i}, &\qquad i&=0,\ldots,\rt-1, &\qquad
   \zt_j&=\frac{\wt_j}{\htilde_j}, &\qquad j&=0,\ldots,\rt-2,
\end{alignat}
\end{subequations}
where
\begin{subequations}\label{def_fgth}
\beq
  f_0=\ldots=f_{r-1}=\gt_0=\ldots=\gt_{\rt-1}=
     \frac{1+w_0}{1+s_0+t_0},
\label{def_fgtha}\eeq
\beq
   h_i=\frac{(1+w_i)(1+w_{i+1})}{(1+s_i+t_i)(1+s_{i+1}+t_{i+1})},
  \qquad i=0,\ldots,r-2\qquad (w_{r-1}\equiv0),
\label{def_fgthb}\eeq
\end{subequations}
\begin{subequations}\label{def_ftght}
\beq
  \ft_0=\ldots=\ft_{\rt-1}=g_0=\ldots=g_{r-1}=
     \frac{1+\wt_0}{1+\st_0+\tt_0},
\label{def_ftghta}\eeq
\beq
   \htilde_i=\frac{(1+\wt_i)(1+\wt_{i+1})}{(1+\st_i+\tt_i)(1+\st_{i+1}+\tt_{i+1})},
  \qquad i=0,\ldots,\rt-2\qquad (\wt_{\rt-1}\equiv0).
\label{def_ftghtb}\eeq
\end{subequations}
\end{prop}

In the notation of Theorem \ref{th:lagrange}, we have 
\beq
\boldsymbol{\xi=(x,y,z,\xt,\yt,\zt)}, \quad
\boldsymbol{\eta=(s,t,w,\st,\tt,\wt)}, \quad 
\boldsymbol{\phi=(f,g,h,\ft,\gt,\htilde)}.
\eeq
\vskip 0.3cm

\noindent
{\bf Proof.}
Consider equation \Ref{eq_ui_xyz} for $u_{i+1}$:
\beq
    u_{i+1}=u_i-x_iu_1u_{i+1}-y_i\ut_1u_{i+1}+z_iu_{i+2}.
\eeq
Note that $z_{i}u_{i+2}=w_{i}u_{i}$
by virtue of \Ref{def_wi}. Having rearranged the terms we get
\beq
    u_{i+1}(1+x_{i}u_1+y_i\ut_1)=u_{i}(1+w_{i}).
\eeq
{}From \Ref{def_s0} it follows that $x_iu_1=s_i$ and
$y_i\ut_1=t_i$. Therefore,
\beq
    \frac{u_i}{u_{i+1}}=\frac{1+s_i+t_i}{1+w_i}.
\label{ui/ui+1}
\eeq
Multiplying the equalities \Ref{ui/ui+1} for the index $i$ and $i+1$
we get
\beq
    \frac{u_i}{u_{i+2}}
   =\frac{(1+s_i+t_i)(1+s_{i+1}+t_{i+1})}{(1+w_i)(1+w_{i+1})}.
\eeq
It remains to replace $u_i/u_{i+2}$ with $z_i/w_i$ from the equality \Ref{def_wi},
and we obtain the equality of the form $z_i=w_i/h_i$, see \Ref{def_fgh}, where
\beq
    h_i=\frac{(1+w_i)(1+w_{i+1})}{(1+s_i+t_i)(1+s_{i+1}+t_{i+1})},
    \qquad i=0,\ldots,r-2,
\label{expr_hi}
\eeq
and it is assumed that $w_{r-1}=0$.

The special case of equation \Ref{eq_ui_xyz} for $i=1$,
\beq
    u_1=1-x_0u_1^2-y_0\ut_1u_1+z_0u_2,
\eeq
has to be treated separately.
Similarly to the general case, we replace $z_0u_2$ with $w_0$ using \Ref{def_w0},
and substitute $u_1=s_0/x_0$ from \Ref{def_s0} and $\ut_1=\st_0/\xt_0$ from
\Ref{def_st0}. Then we replace $y_0\st_0$ with $\xt_0 t_0$,
using both \Ref{def_s0} and \Ref{def_st0}. As a result, the equation takes the form
\beq
    s_0(1+s_0+t_0)=x_0(1+w_0)
\eeq
or, equivalently, the form $x_0=s_0/f_0$, see \Ref{def_fgh} for $i=0$, where
\beq
    f_0=\frac{1+w_0}{1+s_0+t_0}.
\eeq
{}From \Ref{def_s0} and \Ref{def_st0} it also follows that
\beq
    f_i=\frac{s_i}{x_i}= u_1=\frac{\tt_i}{\yt_i}=\gt_i=\frac{s_0}{x_0}=f_0.
\eeq
The remaining half of the equations is obtained
by the tilde-symmetry \Ref{tilde_sym__u}.
\endproof

As was said before, the cases of a corner eigenvalue, i.e. when $r=0$ or $\rt=0$,
and the next one, i.e. when $r=1$ or $\rt=1$, are slightly special. 
When one applies the formulae from the above theorem
describing the generic case, i.e. when both $r,\rt\geq 2$, to such cases
one has to remember that
\begin{itemize}
\item
for a corner eigenvalue, say $\rt=0$, there are no tilded variables and there are
also \textit{no} variables $\boldsymbol{y,t,g}$, so that in this case we only have
$2d-3=2(r+1)-3=2r-1$ variables, therefore one must disregard (\ref{def_ftght})
entirely and set $t_j\equiv0$, $j=0,\ldots,r-1$, in (\ref{def_fgth});
\item
for the immediate next eigenvalue, say $\rt=1$, all variables are present but $\boldsymbol{w}$,
so that one must set $\wt_0\equiv 0$ in (\ref{def_ftghta}) and disregard \Ref{def_ftghtb}
entirely, there are $3d-5=3(r+2)-5=3r+1$ variables in this case.
\end{itemize}
All these modifications are easily seen from the expressions \Ref{def_xyz}--\Ref{def_xyzt} 
of the small parameters $(\boldsymbol{x,y,z,\xt,\yt,\zt})$ in terms of the initial small 
parameters $(\boldsymbol{b,c,\bt,\ct})$. For a more detailed study of the special cases
see section \ref{corner}.

\section{Jacobian}\label{jacobian}
\noindent
In this section only, we use the ordering of variables which is different from the one used above:
$\boldsymbol{\xi}=\boldsymbol{(\xt,\yt,x,y,\zt,z)}$,
$\boldsymbol{\eta=(\st,\tt,s,t,\wt,w)}$,
$\boldsymbol{\phi=(\ft,\gt,f,g,\htilde,h)}$.

\begin{theo}\label{th:jacobian}
The Jacobian $J$ defined by \Ref{defJ} can be expressed as
\beq
J=S\St -T\Tt,
\label{J=SS-TT}
\eeq
where
\beq
S=\frac{1}{W}\left(1+\sum_{j=0}^{r-1}\frac{s_j(1+w_j)}{1+s_j+t_j}\,\prod_{k=0}^{j-1}w_k\right),
\qquad
T=\frac{1}{W}\, \sum_{j=0}^{r-1}\frac{t_j(1+w_j)}{1+s_j+t_j}\,\prod_{k=0}^{j-1}w_k,
\label{sstt}\eeq
and
\beq
W=\prod_{k=0}^{r-2}(1+w_k).
\eeq
\end{theo}

Here we adopt the following agreement: whenever a sum has the upper limit smaller
than the lower one then its value should be taken as 0. Also, any product in a similar
case should be taken as 1.
The tildes in \Ref{J=SS-TT} refer to replacing $s_j$, $t_j$, $w_j$ and $r$ by their
tilded versions. As always, we assume that
$\wt_{\rt-1}\equiv0$ and $w_{r-1}\equiv0$ (see Proposition \ref{pprop}).

Note that the number of terms in $J$ \Ref{J=SS-TT} is
$r\rt+(r+1)(\rt+1)=2r\rt+r+\rt+1=2dk-2k^2+2k-d$ (recall that $d=\rt+r+1$ and $k=\rt+1$).
\vskip 0.3cm

\noindent
{\bf Proof.}
Consider the rows of the matrix
$\mathcal{J}=\delta_{jk}-\frac{\eta_j}{\phi_k}\,\frac{\partial \phi_k}{\partial \eta_j}$:
\beq
\left( 1+\tfrac{\st_0}{1+\st_0+\tt_0},\tfrac{\st_0}{1+\st_0+\tt_0},\ldots,0,\ldots,0,\ldots,
\tfrac{\st_0}{1+\st_0+\tt_0},\ldots,\tfrac{\st_0}{1+\st_0+\tt_0},0,\ldots,0,\ldots\right)
\label{st0}\eeq
\be
\left( 0,1,0,\ldots,0,\ldots,0,\ldots,0,\ldots,
\tfrac{\st_1}{1+\st_1+\tt_1},\tfrac{\st_1}{1+\st_1+\tt_1},0,\ldots,0,\ldots\right)
\ee
\be
\left( 0,0,1,0,\ldots,0,\ldots,0,\ldots,0,\ldots,0,
\tfrac{\st_2}{1+\st_2+\tt_2},\tfrac{\st_2}{1+\st_2+\tt_2},0,\ldots,0,\ldots\right)
\ee
\be
\vdots
\ee
\beq
\left( \tfrac{\tt_0}{1+\st_0+\tt_0},\ldots,1,0,\ldots,0,\ldots,
\tfrac{\tt_0}{1+\st_0+\tt_0},\ldots,\tfrac{\tt_0}{1+\st_0+\tt_0},0,\ldots,0,\ldots\right)
\label{tt0}\eeq
\be
\left( 0,\ldots,0,1,0,\ldots,0,\ldots,0,\ldots,
\tfrac{\tt_1}{1+\st_1+\tt_1},\tfrac{\tt_1}{1+\st_1+\tt_1},0,\ldots,0,\ldots\right)
\ee
\be
\left( 0,\ldots,0,0,1,0,\ldots,0,\ldots,0,\ldots,0,
\tfrac{\tt_2}{1+\st_2+\tt_2},\tfrac{\tt_2}{1+\st_2+\tt_2},0,\ldots,0,\ldots\right)
\ee
\be
\vdots
\ee
\beq
\left( 0,\ldots,\tfrac{s_0}{1+s_0+t_0},\ldots,1+\tfrac{s_0}{1+s_0+t_0},
\tfrac{s_0}{1+s_0+t_0},\ldots,0,\ldots,0,\ldots,\tfrac{s_0}{1+s_0+t_0},0,\ldots\right)
\label{s0}\eeq
\be
\left( 0,\ldots,0,\ldots,0,1,0,\ldots,
0,\ldots,0,\ldots,\tfrac{s_1}{1+s_1+t_1},\tfrac{s_1}{1+s_1+t_1},0,\ldots\right)
\ee
\be
\left( 0,\ldots,0,\ldots,0,0,1,0,\ldots,
0,\ldots,0,\ldots,0,\tfrac{s_2}{1+s_2+t_2},\tfrac{s_2}{1+s_2+t_2},0,\ldots\right)
\ee
\be
\vdots
\ee
\beq
\left( 0,\ldots,\tfrac{t_0}{1+s_0+t_0},\ldots,
\tfrac{t_0}{1+s_0+t_0},\ldots,1,0,\ldots,0,\ldots,\tfrac{t_0}{1+s_0+t_0},0,\ldots\right)
\label{t0}\eeq
\be
\left( 0,\ldots,0,\ldots,0,\ldots,
0,1,0,\ldots,0,\ldots,\tfrac{t_1}{1+s_1+t_1},\tfrac{t_1}{1+s_1+t_1},0,\ldots\right)
\ee
\be
\left( 0,\ldots,0,\ldots,0,\ldots,0,0,1,
0,\ldots,0,\ldots,0,\tfrac{t_2}{1+s_2+t_2},\tfrac{t_2}{1+s_2+t_2},0,\ldots\right)
\ee
\be
\vdots
\ee
\beq
\left(
-\tfrac{\wt_0}{1+\wt_0},\ldots, 0,\ldots,0,\ldots,-\tfrac{\wt_0}{1+\wt_0},\ldots,
\tfrac{1}{1+\wt_0},0,\ldots, 0,\ldots
\right)
\label{wt0}\eeq
\be
\left(
0,\ldots, 0,\ldots,0,\ldots,0,\ldots,-\tfrac{\wt_1}{1+\wt_1},
\tfrac{1}{1+\wt_1},0,\ldots, 0,\ldots
\right)
\ee
\be
\left(
0,\ldots, 0,\ldots,0,\ldots,0,\ldots,0,-\tfrac{\wt_2}{1+\wt_2},
\tfrac{1}{1+\wt_2},0,\ldots, 0,\ldots
\right)
\ee
\be
\vdots
\ee
\beq
\left(
0,\ldots,-\tfrac{w_0}{1+w_0},\ldots, -\tfrac{w_0}{1+w_0},\ldots, 0,\ldots,0,\ldots,
\tfrac{1}{1+w_0},0,\ldots
\right)
\label{w0}\eeq
\be
\left(
0,\ldots, 0,\ldots,0,\ldots,0,\ldots,0,\ldots,-\tfrac{w_1}{1+w_1},
\tfrac{1}{1+w_1},0,\ldots
\right)
\ee
\be
\left(
0,\ldots, 0,\ldots,0,\ldots,0,\ldots,0,\ldots,0,-\tfrac{w_2}{1+w_2},
\tfrac{1}{1+w_2},0,\ldots
\right)
\ee
\be
\vdots
\ee
First lines of each of the four first sub-sets above, which correspond
to the variables $\st_0$ (\ref{st0}), $\tt_0$ (\ref{tt0}), $s_0$ (\ref{s0})
and $t_0$ (\ref{t0}), can be replaced by simpler versions (without
changing the determinant) by adding a multiple of the row associated
with $\wt_0$ (\ref{wt0}) or $w_0$ (\ref{w0}). They become as follows:
\be
\left( 1,0,\ldots,0,\ldots,0,\ldots,
0,\ldots,\tfrac{\st_0(1+\wt_0)}{(1+\st_0+\tt_0)\wt_0},0,\ldots,0,\ldots\right)
\label{st0x}\ee
\be
\left(0,\ldots,1,0,\ldots,0,\ldots,
0,\ldots,\tfrac{\tt_0(1+\wt_0)}{(1+\st_0+\tt_0)\wt_0},0,\ldots,0,\ldots\right)
\label{tt0x}\ee
\be
\left( 0,\ldots,0,\ldots,1,0,\ldots,0,\ldots,0,\ldots,\tfrac{s_0(1+w_0)}{(1+s_0+t_0)w_0},0,\ldots\right)
\label{s0x}\ee
\be
\left( 0,\ldots,0,\ldots,0,\ldots,1,0,\ldots,0,\ldots,\tfrac{t_0(1+w_0)}{(1+s_0+t_0)w_0},0,\ldots\right)
\label{t0x}\ee

After the above replacement, the $(3\rt+3r-2)\times (3\rt+3r-2)$ matrix
$\mathcal{J}=\delta_{jk}-\frac{\eta_j}{\phi_k}\,\frac{\partial \phi_k}{\partial \eta_j}$
acquires the following block-matrix form:
\be
\mathcal{J}=\left(\begin{matrix} \mathbf{1}_{(2\rt+2r)\times(2\rt+2r)}&B\cr C&D
\end{matrix}\right), \qquad
\mbox{\rm whence} \qquad J\equiv \det \mathcal{J}=\det(D-CB),
\ee
where the matrices $B$, $C$, $D$ are
\be
\left(\begin{matrix}
\tfrac{\st_0(1+\wt_0)}{(1+\st_0+\tt_0)\wt_0},&0,&0,&0,\ldots,&0,&0\cr
\tfrac{\st_1}{1+\st_1+\tt_1},&\tfrac{\st_1}{1+\st_1+\tt_1},&0,&0,\ldots,&0,&0\cr
0,&\tfrac{\st_2}{1+\st_2+\tt_2},&\tfrac{\st_2}{1+\st_2+\tt_2},&0,\ldots,&0,&0\cr
\vdots&\vdots&\vdots&\vdots&\vdots&\vdots\cr
0,\ldots,&0,&\tfrac{\st_{\rt-1}}{1+\st_{\rt-1}+\tt_{\rt-1}},&0,\ldots,&0,&0\cr
\tfrac{\tt_0(1+\wt_0)}{(1+\st_0+\tt_0)\wt_0},&0,&0,&0,\ldots,&0,&0\cr
\tfrac{\tt_1}{1+\st_1+\tt_1},&\tfrac{\tt_1}{1+\st_1+\tt_1},&0,&0,\ldots,&0,&0\cr
0,&\tfrac{\tt_2}{1+\st_2+\tt_2},&\tfrac{\tt_2}{1+\st_2+\tt_2},&0,\ldots,&0,&0\cr
\vdots&\vdots&\vdots&\vdots&\vdots&\vdots\cr
0,\ldots,&0,&\tfrac{\tt_{\rt-1}}{1+\st_{\rt-1}+\tt_{\rt-1}},&0,\ldots,&0,&0\cr
0,\ldots,&0,&0,&\tfrac{s_0(1+w_0)}{(1+s_0+t_0)w_0},&0,&0\cr
0,\ldots,&0,&0,&\tfrac{s_1}{1+s_1+t_1},&\tfrac{s_1}{1+s_1+t_1},&0\cr
0,\ldots,&0,&0,&0,&\tfrac{s_2}{1+s_2+t_2},&\tfrac{s_2}{1+s_2+t_2}\cr
\vdots&\vdots&\vdots&\vdots&\vdots&\vdots\cr
0,\ldots,&0,&0,&0,\ldots,&0,&\tfrac{s_{r-1}}{1+s_{r-1}+t_{r-1}}\cr
0,\ldots,&0,&0,&\tfrac{t_0(1+w_0)}{(1+s_0+t_0)w_0},&0,&0\cr
0,\ldots,&0,&0,&\tfrac{t_1}{1+s_1+t_1},&\tfrac{t_1}{1+s_1+t_1},&0\cr
0,\ldots,&0,&0,&0,&\tfrac{t_2}{1+s_2+t_2},&\tfrac{t_2}{1+s_2+t_2}\cr
\vdots&\vdots&\vdots&\vdots&\vdots&\vdots\cr
0,\ldots,&0,&0,&0,\ldots,&0,&\tfrac{t_{r-1}}{1+s_{r-1}+t_{r-1}}
\end{matrix}\right),
\ee
\be
C=\left(\begin{matrix}
-\tfrac{\wt_0}{1+\wt_0},\ldots,&0,\ldots,&0,\ldots,&-\tfrac{\wt_0}{1+\wt_0},\ldots\cr
0,\ldots,&0,\ldots,&0,\ldots,&0,\ldots\cr
\vdots&\vdots&\vdots&\vdots\cr
0,\ldots,&0,\ldots,&0,\ldots,&0,\ldots\cr
0,\ldots,&-\tfrac{w_0}{1+w_0},\ldots,&-\tfrac{w_0}{1+w_0},\ldots,&0,\ldots\cr
0,\ldots,&0,\ldots,&0,\ldots,&0,\ldots\cr
\vdots&\vdots&\vdots&\vdots\cr
0,\ldots,&0,\ldots,&0,\ldots,&0,\ldots
\end{matrix}\right),
\ee
\beq
D=\left(\begin{matrix}
\tfrac{1}{1+\wt_0},&0,&0,\ldots,&0,&0,\ldots,&0\cr
-\tfrac{\wt_1}{1+\wt_1},&\tfrac{1}{1+\wt_1},&0,\ldots,&0,&0,\ldots,&0\cr
\vdots&\vdots&\vdots&\vdots&\vdots&\vdots\cr
0,\ldots,&-\tfrac{\wt_{\rt-2}}{1+\wt_{\rt-2}},&\tfrac{1}{1+\wt_{\rt-2}},&0,&0,\ldots,&0\cr
0,\ldots,&0,&0,&\tfrac{1}{1+w_0},&0,&0\cr
0,\ldots,&0,&0,&-\tfrac{w_1}{1+w_1},&\tfrac{1}{1+w_1},&0\cr
\vdots&\vdots&\vdots&\vdots&\vdots&\vdots\cr
0,\ldots,&0,&0,&0,\ldots,&-\tfrac{w_{r-2}}{1+w_{r-2}},&\tfrac{1}{1+w_{r-2}}
\end{matrix}\right).
\label{DD}\eeq

Therefore, the calculation of the Jacobian has been reduced to finding
an explicit formula for the determinant of a smaller, $(\rt+r-2)\times(\rt+r-2)$,
matrix $D-CB$.
Subtraction of the product $CB$ does only change two rows of the
triangular matrix $D$, those associated with the variables $\wt_0$ and $w_0$.
They now become
\begin{eqnarray}
&&\tfrac{\wt_0}{1+\wt_0}\left(
\tfrac{1}{\wt_0}+\tfrac{\st_0(1+\wt_0)}{(1+\st_0+\tt_0)\wt_0}+\tfrac{\st_1}{1+\st_1+\tt_1},\right.\nonumber\\
&&\qquad\qquad\left.
\tfrac{\st_1}{1+\st_1+\tt_1}+\tfrac{\st_2}{1+\st_2+\tt_2},\ldots,
\tfrac{\st_{\rt-2}}{1+\st_{\rt-2}+\tt_{\rt-2}}
+\tfrac{\st_{\rt-1}}{1+\st_{\rt-1}+\tt_{\rt-1}},\right.\nonumber\\
&&\qquad\left.
\tfrac{t_0(1+w_0)}{(1+s_0+t_0)w_0}+\tfrac{t_1}{1+s_1+t_1},\right.\nonumber\\
&&\qquad\qquad\qquad\left.
\tfrac{t_1}{1+s_1+t_1}+\tfrac{t_2}{1+s_2+t_2},\ldots,
\tfrac{t_{r-2}}{1+s_{r-2}+t_{r-2}}
+\tfrac{t_{r-1}}{1+s_{r-1}+t_{r-1}}
\right)
\nonumber\end{eqnarray}
and
\begin{eqnarray}
&&\tfrac{w_0}{1+w_0}\left(
\tfrac{\tt_0(1+\wt_0)}{(1+\st_0+\tt_0)\wt_0}+\tfrac{\tt_1}{1+\st_1+\tt_1},\right.\nonumber\\
&&\qquad\qquad\left.
\tfrac{\tt_1}{1+\st_1+\tt_1}+\tfrac{\tt_2}{1+\st_2+\tt_2},\ldots,
\tfrac{\tt_{\rt-2}}{1+\st_{\rt-2}+\tt_{\rt-2}}
+\tfrac{\tt_{\rt-1}}{1+\st_{\rt-1}+\tt_{\rt-1}}
,\right.\nonumber\\
&&\qquad\left.
\tfrac{1}{w_0}+\tfrac{s_0(1+w_0)}{(1+s_0+t_0)w_0}+\tfrac{s_1}{1+s_1+t_1},\right.\nonumber\\
&&\qquad\qquad\qquad\left.
\tfrac{s_1}{1+s_1+t_1}+\tfrac{s_2}{1+s_2+t_2},\ldots,
\tfrac{s_{r-2}}{1+s_{r-2}+t_{r-2}}
+\tfrac{s_{r-1}}{1+s_{r-1}+t_{r-1}}
\right),
\nonumber\end{eqnarray}
respectively. The rest of the rows in the matrix $D-CB$ are the same as in the matrix
$D$ above, that is they all have only two non-zero elements: one is on the diagonal
and the other one is to the left of it. By adding multiples of the columns,
without changing the determinant,
we can arrange that all those non-diagonal entries in the selected $(\rt-2)+(r-2)$ rows
vanish, allowing us to compute the determinant by reducing it down
to a $2\times2$ determinant.

Indeed, multiply the last column of the matrix $D-CB$
(look at formula (\ref{DD})) by $w_{r-2}$ and add the result to the
penultimate column. The last row has now got only one (diagonal)
element. Hence, we keep the factor $\frac{1}{1+w_{r-2}}$ and reduce the
determinant to its minor, removing the last column and the last row.
By repeating this process, we shall end up with the expression
\be
J=\frac{\wt_0w_0(\mathcal{S}\tilde{\mathcal{S}}-\mathcal{T}\tilde{\mathcal{T}})}
{(1+\wt_0)\cdots(1+\wt_{\rt-2})(1+w_0)\cdots(1+w_{r-2})}\,,
\ee
where
\begin{eqnarray}
\mathcal{S}&=&\tfrac{1}{w_0}+\tfrac{s_0(1+w_0)}{(1+s_0+t_0)w_0}+\tfrac{s_1(1+w_1)}{1+s_1+t_1}+
\sum_{j=2}^{r-2}\tfrac{s_j(1+w_j)\prod_{k=1}^{j-1}w_k}{1+s_j+t_j}+
\tfrac{s_{r-1}\prod_{k=1}^{r-2}w_k}{1+s_{r-1}+t_{r-1}},\nonumber\\
\tilde{\mathcal{S}}&=&\tfrac{1}{\wt_0}+\tfrac{\st_0(1+\wt_0)}{(1+\st_0+\tt_0)\wt_0}+\tfrac{\st_1(1+\wt_1)}{1+\st_1+\tt_1}+
\sum_{j=2}^{\rt-2}\tfrac{\st_j(1+\wt_j)\prod_{k=1}^{j-1}\wt_k}{1+\st_j+\tt_j}+
\tfrac{\st_{\rt-1}\prod_{k=1}^{\rt-2}\wt_k}{1+\st_{\rt-1}+\tt_{\rt-1}},\nonumber\\
\mathcal{T}&=&\tfrac{t_0(1+w_0)}{(1+s_0+t_0)w_0}+\tfrac{t_1(1+w_1)}{1+s_1+t_1}+
\sum_{j=2}^{r-2}\tfrac{t_j(1+w_j)\prod_{k=1}^{j-1}w_k}{1+s_j+t_j}+
\tfrac{t_{r-1}\prod_{k=1}^{r-2}w_k}{1+s_{r-1}+t_{r-1}},\nonumber\\
\tilde{\mathcal{T}}&=&\tfrac{\tt_0(1+\wt_0)}{(1+\st_0+\tt_0)\wt_0}+\tfrac{\tt_1(1+\wt_1)}{1+\st_1+\tt_1}+
\sum_{j=2}^{\rt-2}\tfrac{\tt_j(1+\wt_j)\prod_{k=1}^{j-1}\wt_k}{1+\st_j+\tt_j}+
\tfrac{\tt_{\rt-1}\prod_{k=1}^{\rt-2}\wt_k}{1+\st_{\rt-1}+\tt_{\rt-1}}.\nonumber
\end{eqnarray}
This is apparently equivalent to the statement of the Theorem.
\endproof

\section{Hypergeometric series}\label{hypergeometric}
\noindent
Now we have all the ingredients
for the right-hand side of the Lagrange formula
\Ref{lagrange_expansion} and can calculate the expansion coefficients.
In this section we return to the initial ordering
$\boldsymbol{\xi=(x,y,z,\xt,\yt,\zt)}$,
$\boldsymbol{\eta=(s,t,w,\st,\tt,\wt)}$
and $\boldsymbol{\phi=(f,g,h,\ft,\gt,\htilde)}$.

Let us choose the function $\chi(\etabf)$ in
\Ref{lagrange_expansion} to be a Laurent monomial
$\chi(\etabf)=\boldsymbol{\eta^{q'}}\equiv
\boldsymbol{s^{m'}t^{n'}w^{p'}\st^{\mt'}\tt^{\nt'}\wt^{\pt'}}$,
$\boldsymbol{q'}\in\Zbbd^D=\Zbbd^{3d-5}$.

Define the step function $\s_{ij}$ as
\beq
    \s_{ij}\equiv \left\{
    \begin{matrix}
        0, &\quad i\le j \\
            1, &\quad i>j
    \end{matrix}\right..
\label{def_stepf}
\eeq
We shall use the binomial, trinomial and quadrinomial
coefficients defined for integer $m, n, p$ as
\begin{subequations}
\beq
    \binom{a}{m}\equiv [x^m](1+x)^a, \qquad
    \binom{a}{m, n}\equiv [x^my^n](1+x+y)^a,
\label{bitrinomial}
\eeq
\beq
    \binom{a}{m, n, p}\equiv [x^my^nz^p](1+x+y+z)^a.
\label{quadrinomial}
\eeq
\end{subequations}
Note that the multinomial coefficients are evaluated as
\begin{subequations}
\beq
   \binom{a}{m}=\frac{a(a-1)\ldots(a-m+1)}{m!}, \qquad
   \binom{a}{m, n}=\frac{a(a-1)\ldots(a-m-n+1)}{m!n!},
\eeq
\beq
   \binom{a}{m, n, p}=\frac{a(a-1)\ldots(a-m-n-p+1)}{m!n!p!}
\eeq
\end{subequations}
for $m\ge0$, $n\ge0$, and $p\ge0$ and vanish if $m<0$, or $n<0$, or $p<0$.

Let $\abs{m}=m_0+\ldots+m_{r-1}$ etc., and set
\beq
    p_{-1}\equiv \abs{m}+\abs{\nt}, \qquad
    \pt_{-1}\equiv \abs{\mt}+\abs{n}.
\label{p-1}
\eeq

We shall also assume that $p_{r-1}=\pt_{\rt-1}=0$.
In the rest of the paper we shall frequently present formulae for the untilded
quantities only, assuming, unless otherwise stated, that the tilded versions are obtained
by tilde-symmetry.

\begin{theo}\label{expansion}
The expansion of the monomial $\boldsymbol{\eta^{q'}}$
in $\boldsymbol{\xi^q}$ is given by
\beq
 \boldsymbol{\eta^{q'}}=\sum_{\boldsymbol{q\geq q'}} \Hcal_{\boldsymbol{q'q}}\,\boldsymbol{\xi^q}
\label{eta_H}
\eeq
where
\begin{eqnarray}
   \Hcal_{\boldsymbol{q'q}}&\equiv&
   [\boldsymbol{\xi^q}]\boldsymbol{\eta^{q'}}\equiv
   [\boldsymbol{x^my^nz^p\xt^{\mt}\yt^{\nt}\zt^{\pt}}]
   \boldsymbol{s^{m'}t^{n'}w^{p'}\st^{\mt'}\tt^{\nt'}\wt^{\pt'}}
   \nonumber\\
  &=& \sum_{i=-1}^{r-1}\sum_{\itilde=-1}^{\rt-1}
  \Fcal^i_{\boldsymbol{q'q}}\tilde{\Fcal}^\itilde_{\boldsymbol{q'q}}
  -\sum_{i=0}^{r-1}\sum_{\itilde=0}^{\rt-1}
  \Gcal^i_{\boldsymbol{q'q}}\tilde{\Gcal}^\itilde_{\boldsymbol{q'q}},
\label{def_H}
\end{eqnarray}
\begin{subequations}\label{def_FG}
\begin{align}
   \Fcal^i_{\boldsymbol{q'q}}
   &=\prod_{j=0}^{r-1}\binom{-p_{j-1}-p_j-\d_{ij}}{m_j-m_j^\prime-\d_{ij},\, n_j-n_j^\prime}
   \binom{p_{j-1}+p_j+\d_{ij}-1}{p_j-p_j^\prime-\s_{ij}},
\label{def_Amnp}\\
  \Gcal^i_{\boldsymbol{q'q}}
  &= \prod_{j=0}^{r-1}\binom{-p_{j-1}-p_j-\d_{ij}}{m_j-m_j^\prime,\, n_j-n_j^\prime-\d_{ij}}
   \binom{p_{j-1}+p_j+\d_{ij}-1}{p_j-p_j^\prime-\s_{ij}}
\label{def_Bmnp}
\end{align}
\end{subequations}
(and respective tilded versions).
\end{theo}
\vskip 0.3cm

\noindent
{\bf Proof.}
Using the definitions \Ref{def_fgth} and \Ref{def_ftght}
we get for $\boldsymbol{q=(m,n,p,\mt,\nt,\pt)}$:
\begin{eqnarray}
\boldsymbol{\phi^q(\eta)}
&\equiv&
 \boldsymbol{f^mg^nh^p\ft^{\mt}\gt^{\nt}\htilde^{\pt}} \nonumber\\
  &=&
   f_0^{m_0}\ldots f_{r-1}^{m_{r-1}}
   g_0^{n_0}\ldots g_{r-1}^{n_{r-1}}
   h_0^{p_0}\ldots h_{r-2}^{p_{r-2}} \nonumber\\
  &&\times
   \ft_0^{\mt_0}\ldots \ft_{\rt-1}^{\mt_{\rt-1}}
   \gt_0^{\nt_0}\ldots \gt_{\rt-1}^{\nt_{\rt-1}}
   \htilde_0^{\pt_0}\ldots \htilde_{\rt-2}^{\pt_{\rt-2}} \nonumber\\
  &=&\left(\frac{1+w_0}{1+s_0+t_0}\right)^{\abs{m}+\abs{\nt}+p_0}\,
 \left(\prod_{j=1}^{r-1}\left(\frac{1+w_j}{1+s_j+t_j}\right)^{p_{j-1}+p_j}\right)
  \nonumber\\
  &&\times\left(\frac{1+\wt_0}{1+\st_0+\tt_0}\right)^{\abs{\mt}+\abs{n}+\pt_0}\,
 \left(\prod_{j=1}^{\rt-1}\left(\frac{1+\wt_j}{1+\st_j+\tt_j}\right)^{\pt_{j-1}+\pt_j}\right),
\nonumber\end{eqnarray}
where we assume $w_{r-1}=\wt_{\rt-1}=p_{r-1}=\pt_{\rt-1}=0$.

Using the shorthand notation \Ref{p-1}
we can write down $\boldsymbol{\phi^q(\eta)}$
in the compact form
\be
   \boldsymbol{\phi^q(\eta)}=
    \left(\prod_{j=0}^{r-1}\left(\frac{1+w_j}{1+s_j+t_j}\right)^{p_{j-1}+p_j}\right)
    \left(\prod_{j=0}^{\rt-1}\left(\frac{1+\wt_j}{1+\st_j+\tt_j}\right)^{\pt_{j-1}+\pt_j}\right).
\ee

Using the step function $\s_{ij}$ \Ref{def_stepf}
we rewrite the expressions \Ref{sstt} for the ingredients of the Jacobian
$J$ in the following  equivalent form:
\be
    T=\sum_{i=0}^{r-1}
    \prod_{j=0}^{r-1}
    \frac{w_j^{\s_{ij}}t_j^{\d_{ij}}}{(1+w_j)^{1-\d_{ij}}(1+s_j+t_j)^{\d_{ij}}},
\ee
\be
    S=\sum_{i=-1}^{r-1}
    \prod_{j=0}^{r-1}
    \frac{w_j^{\s_{ij}}s_j^{\d_{ij}}}{(1+w_j)^{1-\d_{ij}}(1+s_j+t_j)^{\d_{ij}}}.
    \label{Salt}
\ee

Substituting the above expressions for $\chi$, $\boldsymbol{\phi^q}$ and $J$ into
the right-hand side of the Lagrange formula  \Ref{lagrange_expansion}
we get
\beq
    \chi\boldsymbol{\phi^q}J=
    \sum_{i=-1}^{r-1}\sum_{\itilde=-1}^{\rt-1}\Fcal^i\tilde{\Fcal}^\itilde
    -\sum_{i=0}^{r-1}\sum_{\itilde=0}^{\rt-1}\Gcal^i\tilde{\Gcal}^\itilde,
\label{AA-BB}
\eeq
where
\begin{eqnarray}
   \Fcal^i&=& \boldsymbol{s^{m'}t^{n'}w^{p'}} \nonumber\\
        &&\times\prod_{j=0}^{r-1}
        \left(\frac{1+w_j}{1+s_j+t_j}\right)^{p_{j-1}+p_j}
        \left(\frac{s_j}{1+s_j+t_j}\right)^{\d_{ij}}
        \frac{w_j^{\s_{ij}}}{(1+w_j)^{1-\d_{ij}}}
        \nonumber\\
   &=&
   \prod_{j=0}^{r-1}s_j^{m_j^\prime+\d_{ij}}t_j^{n_j^\prime}w_j^{p_j^\prime+\s_{ij}}
    \nonumber\\
    &&\times(1+s_j+t_j)^{-p_{j-1}-p_j-\d_{ij}}
      (1+w_j)^{p_{j-1}+p_j+\d_{ij}-1},
\end{eqnarray}
\beq
   \Gcal^i=
   \prod_{j=0}^{r-1}s_j^{m_j^\prime}t_j^{n_j^\prime+\d_{ij}}w_j^{p_j^\prime+\s_{ij}}
    (1+s_j+t_j)^{-p_{j-1}-p_j-\d_{ij}}
      (1+w_j)^{p_{j-1}+p_j+\d_{ij}-1}
\eeq
and $\tilde{\Fcal}^\itilde$, $\tilde{\Gcal}^\itilde$ are their
tilded versions. In the above formulae it is assumed that
$w_{r-1}^{\s_{r-1,r-1}}=0^0\equiv1$.

It remains to take the coefficient of the resulting expression
for $\chi\boldsymbol{\phi^q}J$ at the monomial $\boldsymbol{\eta^q}$.
Using the binomial and trinomial expansions \Ref{bitrinomial}
we get the expressions \Ref{def_Amnp} and \Ref{def_Bmnp}
for the expansion coefficients
$\Fcal^i_{\boldsymbol{q'q}}\equiv [\boldsymbol{s^mt^nw^p}]\Fcal^i$
and
$\Gcal^i_{\boldsymbol{q'q}}\equiv [\boldsymbol{s^mt^nw^p}]\Gcal^i$.
The final expression \Ref{def_H} for
$[\boldsymbol{\xi^q}]\boldsymbol{\eta^{q'}}$
follows then immediately.
\endproof

As formula \Ref{def_H} shows, the monomial $\boldsymbol{\eta^{q'}}$
expands into $\boldsymbol{\xi^q}$ as a finite sum (the number of the terms equals to
the number of the terms in the Jacobian $J$) of Laurent series.
These Laurent series have the uniform structure and can be in fact expressed in terms
of a single standard series.

Let us introduce the function $\Phi(\boldsymbol{\xi}; \boldsymbol{\mu},\boldsymbol{\tilde{\mu}})$
of $D=3d-5$ complex parameters $\boldsymbol{\xi}$ and depending on two integer vectors
$\boldsymbol{\mu}=(\mu_0,\ldots,\mu_{r-1})\in\Zbbd^r$ and
$\boldsymbol{\tilde{\mu}}=(\tilde{\mu}_0,\ldots,\tilde{\mu}_{\rt-1})\in\Zbbd^{\tilde{r}}$,
altogether $r+\rt=d-1$ integer parameters:
\begin{multline}
 \Phi(\boldsymbol{\xi}; \boldsymbol{\mu},\boldsymbol{\tilde{\mu}})=
 \sum_{\boldsymbol{q}\geq0}
 (-1)^{\abs{m}+\abs{n}+\abs{\mt}+\abs{\nt}}\boldsymbol{\xi^q}
 \\
 \times
 \prod_{j=0}^{r-1}\binom{\mu_j+m_j+n_j+p_{j-1}+p_j-1}{m_j\,,n_j,\,p_j}
 \prod_{j=0}^{\rt-1}\binom{\tilde{\mu}_j+\mt_j+\nt_j+\pt_{j-1}+\pt_j-1}{\mt_j\,,\nt_j,\,\pt_j}.
\label{def_Phi}
\end{multline}

Note that, despite the appearance,
$\Phi(\boldsymbol{\xi}; \boldsymbol{\mu},\boldsymbol{\tilde{\mu}})$
does not factorise into tilde and no-tilde factors
because of $p_{-1}=\abs{m}+\abs{\nt}$ \Ref{p-1}.

For $\mu_j,\tilde{\mu}_j\geq1$ the coefficients of the
series \Ref{def_Phi} can be written down in terms of the Pochhammer symbols
$(a)_m=a(a+1)\ldots(a+m-1)$:
\begin{eqnarray}
 \Phi(\boldsymbol{\xi}; \boldsymbol{\mu},\boldsymbol{\tilde{\mu}})&=&
 \sum_{\boldsymbol{q}\geq0}
  \frac{(-\boldsymbol{x})^{\boldsymbol{m}}(-\boldsymbol{y})^{\boldsymbol{n}}\boldsymbol{z}^{\boldsymbol{p}}%
(-\boldsymbol{\xt})^{\boldsymbol{\mt}}(-\boldsymbol{\yt})^{\boldsymbol{\nt}}\boldsymbol{\zt}^{\boldsymbol{\pt}}}%
{\boldsymbol{m!n!p!\mt!\nt!\pt!}}
   \nonumber\\
 &&\times 
 \prod_{j=0}^{r-1}\frac{(\mu_j)_{m_j+n_j+p_j+p_{j-1}}}{(\mu_j)_{p_{j-1}}}
 \prod_{j=0}^{\rt-1}\frac{(\tilde{\mu}_j)_{\mt_j+\nt_j+\pt_j+\pt_{j-1}}}{(\tilde{\mu}_j)_{\pt_{j-1}}},
\label{Phi_poch}
\end{eqnarray}
where $\boldsymbol{m!}\equiv m_0!\ldots m_{r-1}!$ etc.

Formula \Ref{Phi_poch} characterizes $\Phi$ as a generalised
Appel-Horn-type series: the ratios of adjacent coefficients
of the power series are rational functions of the indices
\cite{Bateman}. Note that $\Phi(0)=1$, and the series
converges in a neighbourhood of $0$ \cite{Bateman}.
The function $\Phi$ can also be viewed as an
$\Acal$-hypergeometric function in the sense
of Gelfand-Kapranov-Zelevinsky \cite{GKZ1989,GKZ1990,GKZ1994}
--- we plan to elaborate on this remark in a subsequent paper.

\begin{theo}\label{Phi_expansion}
The expansion of the monomial $\boldsymbol{\eta^{q'}}$
in $\boldsymbol{\xi^q}$ can be expressed in terms of the function
$\Phi(\boldsymbol{\xi}; \boldsymbol{\mu},\boldsymbol{\tilde{\mu}})$
as follows:
\begin{eqnarray}
 \boldsymbol{\eta^{q'}}&=&
 \boldsymbol{\xi^{q'}}\left(
 \sum_{i=-1}^{r-1}\sum_{\itilde=-1}^{\rt-1}
      \boldsymbol{\xi}^{\ffrak_i+\tilde{\ffrak}_\itilde}
      \Phi\bigl(\boldsymbol{\xi};\boldsymbol{\pi}+\boldsymbol{\nu}^i,
      \boldsymbol{\tilde{\pi}}+\boldsymbol{\tilde{\nu}}^\itilde\,\bigr)\right. \nonumber \\
     &&\left.-\sum_{i=0}^{r-1}\sum_{\itilde=0}^{\rt-1}
     \boldsymbol{\xi}^{\gfrak_i+\tilde{\gfrak}_\itilde}
     \Phi\bigl(\boldsymbol{\xi};\boldsymbol{\pi}+\boldsymbol{\nu}^i,
     \boldsymbol{\tilde{\pi}}+\boldsymbol{\tilde{\nu}}^\itilde\,\bigr)\right),
\label{eta_Phi}
\end{eqnarray}
where
\beq
    \boldsymbol{\xi}^{\ffrak_i}=\prod_{j=0}^{r-1}
    x_j^{\d_{ij}}z_j^{\s_{ij}}=
    \left\{\begin{array}{lcl}
    1, & \quad & i=-1 \\
    x_0, & \quad & i=0 \\
    x_iz_0\ldots z_{i-1}, & \quad & i>0
    \end{array}\right.
\label{def_ffrak}
\eeq  
\beq
    \boldsymbol{\xi}^{\gfrak_i}=\prod_{j=0}^{r-1}
    y_j^{\d_{ij}}z_j^{\s_{ij}}=
    \left\{\begin{array}{lcl}
    1, & \quad & i=-1 \\
    y_0, & \quad & i=0 \\
    y_iz_0\ldots z_{i-1}, & \quad & i>0
    \end{array}\right.
\label{def_gfrak}
\eeq
\beq
    \pi_j=p_{j-1}^\prime+p_j^\prime, \qquad
    \nu^i_j=2\s_{i,j-1},\qquad
    j=0,\ldots,r-1
\label{def_pi}
\eeq
and respectively for the tilded versions. 
In \Ref{def_pi} it is assumed that $p_{-1}^\prime=\abs{m'}+\abs{\nt'}$. 
\end{theo}
\vskip 0.3cm
 
\noindent
{\bf Proof.}
Note that the sum over $\boldsymbol{q\geq q'}$ in \Ref{eta_H}
can be safely replaced with the sum over $\boldsymbol{q}\in\Zbbd^D$,
since the multinomial coefficients in \Ref{def_FG} 
vanish for the negative values of the lower indices and thus 
automatically select the correct limits of summation.
Then from \Ref{eta_H} and \Ref{def_H} we can express 
$\boldsymbol{\eta^{q'}}$ as the finite sum, each term corresponding to a
term in the expansion \Ref{J=SS-TT} of the Jacobian $J$:
\beq
 \boldsymbol{\eta^{q'}}=
 \sum_{i=-1}^{r-1}\sum_{\itilde=-1}^{\rt-1}
 F_{\boldsymbol{q'}}^{i\itilde}(\boldsymbol{\xi})
      -\sum_{i=0}^{r-1}\sum_{\itilde=0}^{\rt-1}
  G_{\boldsymbol{q'}}^{i\itilde}(\boldsymbol{\xi})
\eeq
\beq
    F_{\boldsymbol{q'}}^{i\itilde}(\boldsymbol{\xi})
    =\sum_{\boldsymbol{q}\in\Zbbd^D}
    \Fcal^i_{\boldsymbol{q'q}}\tilde{\Fcal}^\itilde_{\boldsymbol{q'q}}
    \boldsymbol{\xi^q}, \qquad
    G_{\boldsymbol{q'}}^{i\itilde}(\boldsymbol{\xi})
    =\sum_{\boldsymbol{q}\in\Zbbd^D}
    \Gcal^i_{\boldsymbol{q'q}}\tilde{\Fcal}^\itilde_{\boldsymbol{q'q}}
    \boldsymbol{\xi^q}.     
\eeq

Consider the single term $F_{\boldsymbol{q'}}^{i\itilde}(\boldsymbol{\xi})$
of the above sum. From \Ref{def_Amnp} it follows that the summation can be restricted to
the values of the indices for which the multinomial coefficients do not vanish:
$m_j\ge m_j^\prime+\d_{ij}$, $n_j\ge n_j^\prime$, $p_j\ge p_j^\prime+\s_{ij}$,
resp.\ $\boldsymbol{\mt\nt\pt}$. Let us shift the summation indices:
$m_j:=m_j+m_j^\prime+\d_{ij}$, $n_j:=n_j+n_j^\prime$, $p_j:=p_j+p_j^\prime+\s_{ij}$,
etc.\ so that the summation runs over $\boldsymbol{q}\ge0$:
\beq
    F_{\boldsymbol{q'}}^{i\itilde}(\boldsymbol{\xi})
    =\boldsymbol{\xi^{q'}\xi}^{\ffrak_i+\tilde{\ffrak}_\itilde}
     \sum_{\boldsymbol{q}\ge0}
    \Fbbd_{\boldsymbol{q'q}}^i\tilde{\Fbbd}_{\boldsymbol{q'q}}^\itilde
    \boldsymbol{\xi^q},
\eeq
\begin{align}
    \Fbbd_{\boldsymbol{q'q}}^i
    =\prod_{j=0}^{r-1}
    &\binom{-p_{j-1}-p_j-p_{j-1}^\prime-p_j^\prime-2\s_{i,j-1}}{m_j,\,n_j}\notag\\
    \times&\binom{p_{j-1}+p_j+p_{j-1}^\prime+p_j^\prime+2\s_{i,j-1}}{p_j}.
\label{def_Fbbd}
\end{align}

In doing so we have used formulae \Ref{def_Amnp}, \Ref{def_ffrak},
and the identity $\d_{ij}+\s_{ij}=\s_{i,j-1}$. Note that shifting $m_j$
in $p_{-1}=\abs{m}+\abs{\nt}$ produces the term
$\sum_{j=0}^{r-1}\d_{ij}=\s_{i,-1}$ which is consistent with the convention
$p_{-1}^\prime=\abs{m'}+\abs{\nt'}$.

Applying the identity
\beq
    \binom{-a}{m,\,n}\binom{a-1}{p}=(-1)^{m+n}\binom{a+m+n-1}{m,\,n,\,p}
\eeq
to \Ref{def_Fbbd} we get another expression for $\Fbbd_{\boldsymbol{q'q}}^i$:
\beq
        \Fbbd_{\boldsymbol{q'q}}^i
    =(-1)^{\abs{m}+\abs{n}}\prod_{j=0}^{r-1}
    \binom{m_j+n_j+p_{j-1}+p_j+p_{j-1}^\prime+p_j^\prime+2\s_{i,j-1}-1}{m_j,\,n_j,\,p_j}
\eeq

Repeating the same steps for $\tilde{\Fbbd}_{\boldsymbol{q'q}}^\itilde$
and $G_{\boldsymbol{q'}}^{i\itilde}(\boldsymbol{\xi})$
and collecting together the obtained expressions we identify the result as
\Ref{eta_Phi}
\endproof

The expansions we are ultimately interested in are those of the monomials
in the components $v_j$ and $\vt_j$ of the eigenvector $V$ \Ref{jacobi},
or of their rescaled versions: $u_j$ and $\ut_j$ \Ref{vfromu}.

\begin{prop}\label{u_expand}
Let $\boldsymbol{u^k}=(u_1^{k_1},\ldots,u_r^{k_r})$, $\kbf\in\Zbbd^r$, and resp.\ for 
$\boldsymbol{\tilde{u}^{\tilde{k}}}$. 
The expansion of the monomial $\boldsymbol{u^k\tilde{u}^{\tilde{k}}}$ in $\boldsymbol{\xi^q}$
is given by the formula
\begin{eqnarray}
 \boldsymbol{u^k\tilde{u}^{\tilde{k}}}&=&
 \sum_{i=-1}^{r-1}\sum_{\itilde=-1}^{\rt-1}
      \boldsymbol{\xi}^{\ffrak_i+\tilde{\ffrak}_\itilde}
      \Phi(\boldsymbol{\xi};\boldsymbol{\pi}+\boldsymbol{\nu}^i,
      \boldsymbol{\tilde{\pi}}+\boldsymbol{\tilde{\nu}}^\itilde) \nonumber \\
     &&-\sum_{i=0}^{r-1}\sum_{\itilde=0}^{\rt-1}
     \boldsymbol{\xi}^{\gfrak_i+\tilde{\gfrak}_\itilde}
     \Phi(\boldsymbol{\xi};\boldsymbol{\pi}+\boldsymbol{\nu}^i,
     \boldsymbol{\tilde{\pi}}+\boldsymbol{\tilde{\nu}}^\itilde),
\label{u_Phi}
\end{eqnarray}
where $\Phi$, $\boldsymbol{\xi^{\ffrak_i}}$, $\boldsymbol{\xi^{\gfrak_i}}$,
$\boldsymbol{\nu}$ are same as in Theorem \ref{Phi_expansion}, 
and
\beq
    \pi_j=k_{j+1}+\ldots+k_r, \qquad
    j=0,\ldots,r-1,
\label{new_pi} 
\eeq
\end{prop}
\vskip 0.3cm
 
\noindent
{\bf Proof.}
Using formulae \Ref{ufromw} we can express $u_j$ as monomials in $x_0$,
$s_0$, $\zbf$ and $\wbf$. Thus
\beq
    \boldsymbol{u^k}=\boldsymbol{x^{-m'}z^{-p'}s^{m'}w^{p'}},
\label{monomial_u}
\eeq
where
\beq
    m_0^\prime=k_1+k_3+\ldots, \qquad 
    m_1^\prime=\ldots=m_{r-1}^\prime=0,
\eeq
\beq
    p_{2j}^\prime=k_{2j+2}+k_{2j+4}+\ldots, \qquad
    p_{2j+1}^\prime=k_{2j+3}+k_{2j+5}+\ldots
\eeq
(similarly for the tilded variables).
Substituting these values of $p_j^\prime$ into the formula 
$\pi_j=p_{j-1}^\prime+p_j^\prime$
\Ref{def_pi} we get \Ref{new_pi}. Note that 
$p_{-1}^\prime\equiv\abs{m'}=m_0^\prime$.
It remains to apply Theorem \ref{Phi_expansion}
to the monomial \Ref{monomial_u}.
\endproof

To obtain the expansion of the monomial $\boldsymbol{v^k\tilde{v}^{\tilde{k}}}$
in $\boldsymbol{bc\bt\ct}$ one needs to express 
the components $v_i$ of the eigenvector $V$
in terms of $\boldsymbol{u\ut}$ using \Ref{vfromu},
then apply Proposition \ref{u_expand},
and finally, using \Ref{def_xyz} and \Ref{def_xyzt},
to express $\xibf$ in terms of $\boldsymbol{bc\bt\ct}$.
To get the expansions involving the eigenvalue
$\la$ one has to use \Ref{lav1}.
In particular, the resulting expansions for $u_i$ or $\ut_i$ give an 
explicit solution for iterations of equations \Ref{eq_uuti_xyz}.

In this paper we study only formal power series expansions.
The analytic version of the Lagrange theorem
\cite{Good1960} guarantees that our power series converge
in an open neighbourhood around $\xibf=\boldsymbol{0}$.
The precise description of the convergence domain
of a multivariate power series is, however, a difficult task,
see for example \cite{GKZ1994,Bateman}, and we leave it to
further study.

\section{Corner eigenvalue}\label{corner}
\noindent
As mentioned in the end of section \ref{lagrange}, the cases 
of the corner eigenvalue $r=0$ ($\rt=0$) have some peculiarities.
Let us consider the case $\rt=0$ in more detail
(the case $r=0$ can be treated by the tilde-symmetry).
We have
\beq
  \begin{pmatrix}
     0 & b_0  \\
     c_0 & a_1 & b_1 \\
    & & \hdotsfor{3} \\
     & & & c_{i-1} & a_i & b_i \\
    & & & & \hdotsfor{3} \\
     & & & & & c_{r-1} & a_r
  \end{pmatrix}
  \begin{pmatrix}
   1 \\ v_1 \\ \hdots \\ v_i \\ \hdots \\ v_r
  \end{pmatrix}
  =\lambda
  \begin{pmatrix}
   1 \\ v_1 \\ \hdots \\ v_i \\ \hdots \\ v_r
  \end{pmatrix}.
\label{jacobi_corner}
\eeq
The tilded variables are absent. Besides, the sequences $y_i$, resp.\
$t_i$ are also absent.
The total set of variables $\boldsymbol{\xi=(x,z)}$
has the cardinality $(r)+(r-1)=2r-1=2d-3$, and all
the expansion variables are
independent, in contrast with the generic case.
Respectively, $\boldsymbol{\eta=(s,w)}$ and
$\boldsymbol{\phi=(f,h)}$, where
\beq
   f_0=\ldots=f_{r-1}=\frac{1+w_0}{1+s_0}, \qquad
   h_i=\frac{(1+w_i)(1+w_{i+1})}{(1+s_i)(1+s_{i+1})},
  \qquad i=0,\ldots,r-2.
\eeq

The expression \Ref{J=SS-TT} for the jacobian $J$ simplifies to
\beq
J=S=\left(1+\sum_{j=0}^{r-1}\frac{s_j(1+w_j)}{1+s_j}\,\prod_{k=0}^{j-1}w_k\right)
    \prod_{k=0}^{r-2}(1+w_k)^{-1}
    \eeq
and contains thus only $r+1=d$ terms. The analog of formula
\Ref{AA-BB} is
\beq
        \chi\boldsymbol{\phi^q}J=
    \sum_{i=-1}^{r-1}\Fcal^i,
\eeq
where
\beq
    \Fcal^i=
    \prod_{j=0}^{r-1}s_j^{m_j^\prime+\d_{ij}}w_j^{p_j^\prime+\s_{ij}}
     (1+s_j)^{-p_{j-1}-p_j-\d_{ij}}
      (1+w_j)^{p_{j-1}+p_j+\d_{ij}-1}.
\eeq

The expression \Ref{def_H}
for the expansion coefficient becomes
\beq
   \Hcal_{\boldsymbol{q'q}}\equiv
   [\boldsymbol{\xi^q}]\boldsymbol{\eta^{q'}}\equiv
   [\boldsymbol{x^mz^p}]
   \boldsymbol{s^{m'}w^{p'}}
   = \sum_{i=-1}^{r-1}
   \Fcal^i_{\boldsymbol{q'q}},
\eeq
where
\beq
   \Fcal^i_{\boldsymbol{q'q}}
   \equiv [\boldsymbol{s^mw^p}]\Fcal^i
   =\prod_{j=0}^{r-1}\binom{-p_{j-1}-p_j-\d_{ij}}{m_j-m_j^\prime-\d_{ij}}
   \prod_{j=0}^{r-2}\binom{p_{j-1}+p_j+\d_{ij}-1}{p_j-p_j^\prime-\s_{ij}},
\eeq
whith $p_{-1}=\abs{m}$.

The expansion of the monomial $\boldsymbol{u^k}$ in $\boldsymbol{x^mz^p}$
is given by the analog of formula \Ref{u_Phi}:
\beq
 \boldsymbol{u^k}=
 \sum_{i=-1}^{r-1}
      \boldsymbol{\xi}^{\ffrak_i}
      \Phi(\boldsymbol{x,z};\boldsymbol{\pi}+\boldsymbol{\nu}^i)
\label{u_Phi_corner}
\eeq
where $\boldsymbol{\xi}^{\ffrak_i}$, $\boldsymbol{\nu}^i$
are same as in Theorem \ref{Phi_expansion},
$\boldsymbol{\pi}$ is given by \Ref{new_pi}, 
and $\Phi(\boldsymbol{x,z};\boldsymbol{\mu})$ is defined as the series
\beq
 \Phi(\boldsymbol{x,z};\boldsymbol{\mu})=
 \sum_{\boldsymbol{m,p}\geq0}
 \boldsymbol{(-x)^mz^p}
 \prod_{j=0}^{r-1}\binom{\mu_j+m_j+p_{j-1}+p_j-1}{m_j,\,p_j},
\label{def_Phi_corner}
\eeq
or, for $\mu_j\ge1$,
\beq
 \Phi(\boldsymbol{x,z};\boldsymbol{\mu})=
 \sum_{\boldsymbol{m,p}\geq0}
  \frac{(-\boldsymbol{x})^{\boldsymbol{m}}\boldsymbol{z}^{\boldsymbol{p}}}{\boldsymbol{m!p!}}
 \prod_{j=0}^{r-1}\frac{(\mu_j)_{m_j+p_j+p_{j-1}}}{(\mu_j)_{p_{j-1}}}.
\label{Phi_poch_corner}
\eeq

\section{Examples}\label{examples}
\noindent
In this section we illustrate our general results with a few low-dimensional examples.

{\bf 6.1} {\em The simplest case is} $d=2$ ($r=1$, $\rt=0$). The spectral problem is
\beq
    \begin{pmatrix}0 & b_0 \\ c_0 & a_1\end{pmatrix}
    \begin{pmatrix}1 \\ v_1\end{pmatrix}
    =\la\begin{pmatrix}1 \\ v_1\end{pmatrix},
\eeq

After eliminating the eigenvalue $\la=b_0v_1$ and introducing the rescaled variables
$u\equiv u_1=-v_1a_1/c_0$, $x\equiv x_0=b_0c_0/a_1^2$ we get for the single unknown
variable $u$ the single quadratic equation
\beq
    1-u-xu^2=0.
\eeq

The branch we are studying is selected by the condition $u\rightarrow1$ for $x\rightarrow0$,
and the corresponding solution is
\beq
    u(x)=\frac{\sqrt{1+4x}-1}{2x}.
\eeq

In terms of the variable $s\equiv xu$ we get the equivalent Lagrange-type equation
\beq
    x=s(1+s)=\frac{s}{f}, \qquad f=\frac{1}{1+s}.
\eeq

The corresponding Jacobian $J$ has two terms:
\beq
    J=1-\frac{s}{f}\,\frac{df}{ds}=1+\frac{s}{1+s}.
\eeq

The expansion of the monomial $u^k$, $k\in\Zbbd$ given by Proposition \ref{u_expand} is
\beq
    u^k=\Phi(x;k)+x\Phi(x;k+2),
\label{d2uk}
\eeq
where
\beq
    \Phi(x;\mu)=\sum_{m\ge0} (-x)^m\,\binom{\mu+2m-1}{m}
    =\frac{1}{\sqrt{1+4x}}\left(\frac{\sqrt{1+4x}-1}{2x}\right)^{\mu-1},
\eeq

{}For $\mu\ge1$ the function $\Phi(x;\mu)$ can be expressed in terms of the Gauss hypergeometric
series:
\beq
    \Phi(x;\mu)=\sum_{m\ge0} \frac{(-x)^m}{m!}\,\frac{(\mu)_{2m}}{(\mu)_m}
    ={}_2F_1\left(\begin{array}{c}\frac{\mu}{2},\,\frac{\mu+1}{2} \\ \mu \end{array}; -4x\right).
\eeq

As a matter of fact, for $k\geq0$ the two hypergeometric series in \Ref{d2uk} can be summed into a single series
\beq
    u^k=\sum_{m\ge0} \frac{(-x)^m}{m!}\,\frac{(k)_{2m}}{(k+1)_m}
    ={}_2F_1\left(\begin{array}{c}\frac{k}{2},\,\frac{k+1}{2} \\ k+1 \end{array}; -4x\right)
    =\left(\frac{\sqrt{1+4x}-1}{2x}\right)^k.
\eeq

\vskip1mm
{\bf 6.2} {\em Three-dimensional case, corner eigenvalue:} $d=3$ ($r=2$, $\rt=0$).
The spectral problem
\beq
    \begin{pmatrix}
        0 & b_0 & 0\\
        c_0 & a_1 & b_1 \\
        0 & c_1 & a_2
    \end{pmatrix}
    \begin{pmatrix}
        1 \\ v_1 \\ v_2
    \end{pmatrix}
    =\la\begin{pmatrix}1 \\ v_1 \\ v_2\end{pmatrix},
\eeq
after the substitutions $\la=b_0v_1$, $v_1=-u_1c_0/a_1$, $v_2=u_2c_0c_1/a_1a_2$
is reduced to the pair of quadratic equations for the variables $u_1$, $u_2$
\begin{subequations}
\begin{eqnarray}
   1-u_1-x_0u_1^2+z_0u_2&=&0,\\
   u_1-u_2-x_1u_1u_2&=&0,
\end{eqnarray}
\end{subequations}
where the expansion variables are
\beq
   x_0=\frac{b_0c_0}{a_1^2}, \qquad
   x_1=\frac{b_0c_0}{a_1a_2}, \qquad
   z_0=\frac{b_1c_1}{a_1a_2},
\eeq
and we choose the branch $u_1,u_2\rightarrow1$ as $x_0,x_1,z_0\rightarrow0$.

The equivalent Lagrange-type equations for the variables
$s_0=x_0u_1$, $s_1=x_1u_1$, $w_0=z_0u_2$ are
\beq
   x_0=\frac{s_0}{f_0}, \qquad
   x_1=\frac{s_1}{f_1}, \qquad
   z_0=\frac{w_0}{h_0},
\eeq
where
\beq
f_0=f_1=\frac{1+w_0}{1+s_0}, \qquad
h_0=\frac{1+w_0}{(1+s_0)(1+s_1)}.
\eeq

The three-term Jacobian $J$ equals
\beq
    J=\frac{1}{1+w_0}+\frac{s_0}{1+s_0}+\frac{s_1w_0}{(1+s_1)(1+w_0)}.
\eeq

By formula \Ref{u_Phi_corner},  
the monomial $u_1^{k_1}u_2^{k_2}$ expands in $x_0$, $x_1$, $z_0$ as
\begin{align}
    u_1^{k_1}u_2^{k_2}&=\Phi(x_0,x_1,z_0;k_1+k_2,k_2)\notag\\
    &+x_0\Phi(x_0,x_1,z_0;k_1+k_2+2,k_2)\notag\\
    &+x_1z_0\Phi(x_0,x_1,z_0;k_1+k_2+2,k_2+2),
\end{align}
where
\begin{multline}
    \Phi(x_0,x_1,z_0;\mu_0,\mu_1)=\sum_{m_0,m_1,p_0\ge0}
    (-x_0)^{m_0}(-x_1)^{m_1}z_0^{p_1}\\
    \times\binom{\mu_0+2m_0+m_1+p_0-1}{m_0,\,p_0}
    \binom{\mu_1+m_1+p_0-1}{m_1},
\end{multline}
or, for $\mu_0,\mu_1\ge1$,
\beq
  \Phi(x_0,x_1,z_0;\mu_0,\mu_1)=\sum_{m_0,m_1,p_0\ge0}
  \frac{(-x_0)^{m_0}(-x_1)^{m_1}z_0^{p_1}}{m_0!m_1!p_0!}\,
  \frac{(\mu_0)_{2m_0+m_1+p_0}(\mu_1)_{m_1+p_0}}{(\mu_0)_{m_0+m_1}(\mu_1)_{p_0}}.
\eeq

\vskip1mm
{\bf 6.3} {\em Three-dimensional case, middle eigenvalue:} $d=3$ ($r=\rt=1$).
The spectral problem
\beq
    \begin{pmatrix}
      \at_1 & \bt_0 & 0 \\
      \ct_0 & 0 & b_0 \\
      0 & c_0 & a_1
    \end{pmatrix}\,
    \begin{pmatrix}
       \vt_1 \\ 1 \\ v_1
    \end{pmatrix}
    =\la\,
    \begin{pmatrix}
       \vt_1 \\ 1 \\ v_1
    \end{pmatrix}
\eeq
after the substitutions 
$\la=\ct_0\vt_1+b_0v_1$, $v_1=-uc_0/a_1$, $\vt_1=-\ut\bt_0/\at_1$
is reduced to the pair of quadratic equations for the variables $u$, $\ut$
\begin{subequations}
\begin{eqnarray}
     1-u-xu^2-y\ut u&=&0, \\
     1-\ut-\xt\ut^2-\yt u\ut&=&0,
\end{eqnarray}
\end{subequations}
where the expansion variables are
\beq
     x=\frac{b_0c_0}{a_1^2}, \qquad
     y=\frac{\bt_0\ct_0}{\at_1a_1}, \qquad
     \xt=\frac{\bt_0\ct_0}{\at_1^2}, \qquad
     \yt=\frac{b_0c_0}{a_1\at_1},
\eeq
and the chosen branch is $u,\ut\rightarrow1$ as $x,y,\xt,\yt\rightarrow0$.
Note that the variables $x,y,\xt,\yt$ are bound by the single relation
$x\xt=y\yt$.

The equivalent Lagrange-type equations for the variables
$s=xu$, $\st=\xt\ut$, $t=y\ut$, $\tt=\yt u$ are
\beq
   x=\frac{s}{f}\,, \qquad
   y=\frac{t}{g}\,, \qquad
   \xt=\frac{\st}{\ft}\,, \qquad
   \yt=\frac{\tt}{\gt}\,,
\eeq
where
\beq
  f=\gt=\frac{1}{1+s+t}\,, \qquad
  \ft=g=\frac{1}{1+\st+\tt\,}\,.
\eeq

The five-term Jacobian $J$ equals
\begin{multline}
    J=1+\frac{s}{1+s+t}+\frac{\st}{1+\st+\tt\,}\\
    +\frac{s\st}{(1+s+t)(1+\st+\tt\,)}
    -\frac{t\tt}{(1+s+t)(1+\st+\tt\,)}.
\end{multline}

By formula \Ref{u_Phi},  
the monomial $u^k\ut^{\kt}$ expands in $\boldsymbol{\xi}=(x,y,\xt,\yt)$ as
\begin{align}
    u^k\ut^{\kt}&=\Phi(\boldsymbol{\xi};k,\kt)
       +x\Phi(\boldsymbol{\xi};k+2,\kt)
       +\xt\Phi(\boldsymbol{\xi};k,\kt+2)\\
       &+x\xt\Phi(\boldsymbol{\xi};k+2,\kt+2)
       -y\yt\Phi(\boldsymbol{\xi};k+2,\kt+2),
\end{align}
where, for $\qbf=(m,n,\mt,\nt)$,
\begin{multline}
   \Phi(\boldsymbol{\xi};\mu,\tilde{\mu})=
     \sum_{\qbf\ge0}(-x)^m(-y)^n(-\xt)^{\mt}(-\yt)^{\nt} \\
     \times\binom{\mu+2m+n+\nt-1}{m,\,n}
     \binom{\tilde{\mu}+2\mt+n+\nt-1}{\mt,\,\nt},
\end{multline}
or, for $\mu,\tilde{\mu}\ge1$,
\beq
   \Phi(\boldsymbol{\xi};\mu,\tilde{\mu})=
     \sum_{\qbf\ge0}
     \frac{(-x)^m(-y)^n(-\xt)^{\mt}(-\yt)^{\nt}}{m!n!\mt!\nt!}\,
     \frac{(\mu)_{2m+n+\nt}(\tilde{\mu})_{2\mt+n+\nt}}{(\mu)_{m+\nt}(\tilde{\mu})_{\mt+n}}.
\eeq

\section{Discussion}\label{discussion}
\noindent
In the present paper we have solved a problem of much physical
interest. The spectra of finite Jacobi matrices appear in many applications:
from orthogonal polynomials and nearest-neighbours interaction models
to solvable models of quantum mechanics (Lam\'e polynomials and
Bethe ansatz). 

What is left for further study are the questions of convergence domains, differential
equations and integral representations for the obtained hypergeometric series, as well as
their relation to $\Acal$-hypergeometric functions introduced in \cite{GKZ1989,GKZ1990,GKZ1994}.
The approach used in our work can be generalised to multiparameter spectral
problems \cite{Sleeman}.

Explicit perturbative solutions have also appeared in a different context in the works by
Edwin Langmann \cite{Langmann1,Langmann2}, notably for the multi-dimensional spectral problems
related to the Calogero-Sutherland and elliptic Calogero-Moser
systems. 

\section*{Acknowledgements}
\noindent

This work has been partially supported by the European Community (or European Union) through the FP6 
Marie Curie RTN {\sl ENIGMA} (Contract number MRTN-CT-2004-5652).


\end{document}